\documentstyle[12pt]{article}

\begin{document} 
\begin{center} {\large\bf
Notes to the early history of the Knot Theory in Japan.}
\end{center}
\begin{center}
To be published in Annals of the Institute for 
Comparative Studies of Culture; Tokyo Woman's Christian University,
under the title:\ The interrelation of the Development of Mathematical 
Topology in Japan, Poland and USA.
\end{center}
\begin{center}
J\'ozef H. Przytycki 
\end{center}

The goal of this essay is to give a preliminary description
of the growth of research on Knot Theory 
in Japan and to look at its origins. In particular, the influences of 
R.H. Fox of Princeton is analyzed\footnote{Ralph Hartzler Fox (1913-1973).
We give his short biography in Section 7.}.

\section{Early history of Topology in Japan}
The University of Tokyo\footnote{In the period 1886-1947, 
Tokyo Imperial University. We say concisely Tokyo University.}
and the Tokyo Mathematical Society 
were founded in 1877, 
the tenth year of the Meiji era (1868-1912).
During the period 1897-1942 six more Imperial 
Universities were founded: Kyoto 1897, 
T\^{o}hoku (Sendai) 1907\footnote{By 1911, 
the Tohoku Mathematical Journal was initiated. It was Japan's first 
western style, research-level journal.}, 
Hokkaido (Sapporo) 1930, Osaka 1931, Kyushu (Fukuoka) 1939, and Nagoya 1942. 
In addition, two other universities (called Higher 
Normal Schools\footnote{They were 
professional training schools with the main goal to produce middle 
school teachers.}), later
renamed Bunrika Daigaku, with 
departments of mathematics were founded in Tokyo and Hiroshima in 1929.

   Most likely Takeo Wada (1882-1944) published the first 
paper in Japan devoted to topology \cite{Wad}, 1911/1912.
Wada graduated from Kyoto University and became an assistant professor 
in 1908.  He visited the USA, France and Germany from 1917 to 1920.

Kuniz\^o Yoneyama (1877-1968) did research in topology under Wada's 
influence \cite{Wad}.
He wrote (1917-1920) the first Japanese tract on General Topology.

Wada's and Yoneyama's work was continued by Hidetaka Terasaka (see Section 2).

In 1929 K\^oshiro Nakamura (1901-1985) from Tokyo Bunrika University 
(Higher Normal School),
visited Berlin University and studied algebraic topology under H.Hopf.
Then, following Hopf, he visited Zurich (Switzerland). He came back to Japan
in 1932 and published two books introducing algebraic topology to Japan.
Atuo Komatu (1909-1995) was then a third year student at University 
of Tokyo, and he was strongly influenced by Nakamura.
Komatu graduated from Tokyo University in 1932. He wrote his 
doctoral thesis under the supervision of S.Nakagawa (see the 
footnote 11) in 1941, \cite{Kom}.

Work on algebraic topology was also done by Makoto Abe (1914-1945) and
Kunihiko Kodaira (1915-1997), both of whom were students of  
Iyanaga\footnote{Sh\^{o}kichi Iyanaga was born in 1906 and
graduated from Tokyo University in 1929 (he studied with T.Takagi). 
He spent 3 years in Europe 
(1931-34), partially in Germany (Hamburg) where he was a student of E.Artin. 
In 1937, after his teacher, Nakagawa, reached 60 years of age, 
the {\it Chair} in geometry was passed to Iyanaga.
Iyanaga was a specialist 
in number theory and geometry but was also
interested in topology, for example he gave a talk {\it The foundation of
general Topology} at the 4th Topology Colloquium at Kyoto University in
April 1939. He was also a member of the editorial committee of 
the bulletin of the Topology Colloquium (from April 1938) \cite{Mim};
compare the entry on Murasugi.}
 at Tokyo University. Abe died in the chaos after the end of
World War II\footnote{He was of week health and was overworked 
when in a charge of an evacuation of Tokyo University. Lung cancer and
malnutrition caused his death in 1945.}. Kodaira received the Fields medal in 1954.

Kinjir\^o Kunugui (1903-1975) went in 1928 to study in Strasbourg, 
and then in Paris under Fr\'echet 
and he did his PhD in Paris on Dimension Theory.

The fundamental book on Dimension Theory by Hurewicz and Wallman (1941) 
motivated research of Kiiti Morita (1915--1995), Jun-iti Nagata\\
 (1925--), and Kei\^o Nagami (1925--).

Shizuo Kakutani (1911--) graduated in 1934 from T\^{o}hoku 
University and became an
assistant at the newly founded mathematics department of Osaka University,
where he was under the general guidance of T.Shimura. 
His 1936 paper on Riemannian 
surfaces caught the attention of H.Weyl who invited Kakutani 
to IAS (Princeton) in 1940. 
The war forced him to go back to Japan in 1942.

In next sections we describe people involved in the early study of knot theory
in Japan. We follow, in our description, the ``generation" approach
of Noguchi \cite{No-3}. This period is summarized by M.Mimura as follows
\cite{Mim}: ``Right after World War II, Terasaka and Kinoshita in Osaka
and Noguchi, Homma and Murasugi in Tokyo were working by themselves without
good contacts, but their research made good progress when 
Homma and Noguchi paid a visit to Osaka University. 
Around 1950 Terasaka, Kinoshita and Homma 
were studying homeomorphisms of $R^2$ and of $R^3$... 
In this research the importance of how 
to embed a closed curve (knot) and a surface (2-dimensional manifold)
in $R^3$ attracted attention.\footnote{Another reason for interest in knots by
Kinoshita and Terasaka may be their interest in transformations of 
the Euclidean 3-space and their fixed point sets.}
 Homma proved Dehn's lemma \cite{Hom-2}), i.e. whether 
a knot ($S^1 \subset R^3$) is truly unknotted is determined by whether 
the knot group $\pi_1(R^3 - S^1$) is isomorphic to $Z$, 
almost at the same time as,
and independently of, Papakyriakopoulos. Terasaka, Kinoshita 
and others also undertook research in knot theory. One of the typical results 
at that time was the generalization of composition of knots by 
Kinoshita and Terasaka \cite{K-T}. 
It can be said that research in knot theory in Japan
truly started around this time".
Several influential
paper on knot theory were published in Japan in 1957 -1958 period
\cite{Has,Ha-Ho,Hom-2,Hos,Kin-4,Kin-5,Kin-6,Kin-7,Kin-8,K-T,Mur-1,
Mur-2,Saka,Ya-Ki,Fo-1,Fo-2}. The earliest of these papers \cite{Kin-4}
was received by Osaka Mathematical Journal on March 25, 
1957\footnote{Most likely the first in Japan course on knot theory 
(based on Reidemeister's book) was conducted by Gaisi Takeuti at
Tokyo University of Education in 1952 or 1953. K.Murasugi 
attended this course, compare Section 3.}.
 Mimura adds: ``Visits to Japan by people like 
Moise, Fox, Cairns, Bing, Harrold and visits to the USA by Noguchi,
Kinoshita, Homma, Murasugi, Hosokawa, Junzo Tao made possible a close
research cooperation and the exchange of information between Japan and the
USA"\footnote{One should add here frequent visits by J.P.Mayberry
\cite{Mur-4}; compare the entry on Murasugi.}. 
Visits by Fox to Japan\footnote{His first visit
to Japan was as an American representative to the World {\it Go}-game 
tournament in October 1963. The universities of Tokyo, Osaka, Nagoya, 
Kyoto and Fukuoka have invited him to lecture on knot theory
 while he was in Japan \cite{Fo-9}; compare Section 7
.}
 and Japanese knot theoretists
to Princeton (where Fox worked) were of great importance. Fox was
very friendly and encouraging toward young topologists and often 
invited them to Princeton. 

\section{First generation.}

Hidetaka Terasaka, born 1904 in Tokyo, died April 3, 1996.\\
He graduated from Tokyo University in 1928, as a student
of Nakagawa\footnote{Senkichi Nakagawa, (1876-1942) was
born in Kanazawa and graduated from
the mathematics department of Tokyo University (1898).
He went to Germany where he studied geometry. After returning to Japan,
he was appointed assistant professor at his alma mater, and professor (1914).
He was a specialist in Non-Euclidean geometry \cite{Ja}.}.
In 1933 Terasaka went to Germany and Austria (Vienna)
as a scholarship fellow of the Ministry of Education
of the Japanese Government. In Vienna, where he spent more than
a year, he studied mostly with Karl Menger (1902-1985)\footnote{
Menger held the chair of geometry at the University of Vienna (after
Reidemeister left for K\"{o}nigsberg) from 1927 till 1938 \cite{Sig}.}. 
He met there also 
Heinrich Tietze (1880-1964), and Herbert Seifert (1906-1996), 
who just came from Dresden\footnote{
Noguchi wrote in his stories of mathematicians \cite{No-3}: 
``In the 20'th there were 
two trends in topology, one represented by Brouwer, Uryson and Menger, 
second by P.Heegaard, S.B.Ker\'ekj\'arto and Veblen. 
Seifert followed the second trend."}.
Terasaka was supposed to stay in Vienna for 2 years but came back after
one year to become an assistant professor (and a year later a full professor)
in the newly opened Osaka University. 
Other full professors, in addition to Terasaka, at Osaka in 1935 were:
Kunugi, Shimizu, Shoda and K\^osaku Yoshida. A few years later they were 
joined by Kakutani, Komatu, Nakamura and Sakata (Shizuma\footnote{
Ryozi Sakata (1916-1985) was adopted by the Shizuma family in 1942,
after which he used the name Shizuma \cite{Mim}.}).
Terasaka was appointed to the editorial committee of the Bulletin
of the Topology Colloquium (April 1938), which was then named Isosugaku
(Japanese name for Topology). The first issue was
published in October 1938 \cite{Mim}.

Terasaka was familiar with the work of the Polish school of topology and
in particular he applied to general lattices Kuratowski's method of defining
a topological space by using the closure operator. He submitted his paper
to the leading Polish mathematical journal 
{\it Fundamenta Mathematicae}\footnote{Fundamenta Mathematicae was created
by young Polish Mathematicians, Z.Janiszewski (1888-1920), 
S.Mazurkiewicz (1888-1945) and W.Sierpi'nski (1882-1969), 
in 1920, after Poland gained independence in 1918. In a sense, 
the situation of Poland and Japan was similar at that time with respect 
to grow of mathematics: several young mathematicians, 
educated in Germany or France were returning to their
native countries and building, almost from scratch, schools of
modern mathematics.}.
It was accepted for publication just before the Second World War.
However Poland was attacked by nazis on the first day of September of 
1939 and the journal was burned\footnote{
K.Kuratowski writes in \cite{Kur}:``Before the war `Fundamenta' published 
32 volumes. Volume XXXIII was in print. The first several
pages were printed before the war.
Most of the materials and manuscripts were destroyed on order of
the director of the printing shop [in old Jagiellonian University
in Cracow], who had been nominated to that position by the nazis.  
However, due to the exceptional efforts of workers of the printing shop,
a lot of materials (proofs, manuscripts) were hidden and survived till
the end of the war...
Volume XXXIII of 'Fundamenta', containing to great extend, recovered prewar
results and some results obtained in the war period, was published before 
the end of 1945."}. Fortunately off-prints (50 copies) safely reached 
Terasaka in Japan.  In the Fundamenta Mathematicae  
issue (with date 1939) published in Poland after 
the war only the title of the Terasaka paper is included as he thought 
that after six years his paper was not up to date, and reprints were already 
sent out \cite{Ter-1}\footnote{One can add here that Kinoshita, Terasaka's
student, also published his early papers in {\it Fundamenta} 
\cite {Kin-2,Kin-3}.}.

Terasaka heard talks on the knot theory by Seifert in Vienna, but he
started to work on the topic by himself later, after 1953.
Most likely in this he was influenced by his student, Shin'ichi Kinoshita 
(see the entry on Kinoshita).
Some other Terasaka students were Tao and Hosokawa (Section 3).
In 1961 Terasaka moved from Osaka to Tokyo Women
Christian University (TWCU) invited by Aiko Negishi his former student 
(she graduated from TWCU but went 
to Osaka University to study with Terasaka). 
After 2-3 years he moved to Sophia
University to the newly founded mathematics department, invited 
by his old colleague from Tokyo University, M.Moria (later 
president of Sophia University). Hisako Kondo was his master 
degree student at Sophia University (1975)\footnote{She was born 
in 1951 and graduated in March of 1973 from
Japan Women's University (Nihon Jyosidaigaku). She wrote an influential
paper on the Alexander polynomial of the unknotting number one 
knots \cite{Kon}}.
Junzo Wada\footnote{He was a professor at TWCU and later at Waseda
University. Now he is retired.}
 was another student of Terasaka specializing in functional analysis.
Takaaki Yanagawa and Takeshi Yajima were also students of Terasaka 
at Osaka University.

Terasaka died April 3, 1996.

\section{Second generation}

From the Preface to \cite{Su-2}: ``The Topology and Computer 
Science Symposium was held at Akanezaki Hotel in Atami 
from the 17th to the 20th January, 1986.
The Symposium was in honour of Professors Shin'ichi Kinoshita, 
Hiroshi Noguchi and Tatsuo Homma on their sixtieth birthday.
The three professors studied topology at the Osaka University under
Professor Hidetaka Terasaka. During the past 35 years their ardor
and skills have contributed to Knot Theory, PL Topology and 
Computer Science, and they educated many students. Most of the
participants of the symposium were their former students who
owe much to their encouragement and stimulus."

\begin{enumerate}
\item[(K)] Shin'ichi Kinoshita, born 1925.\\
He was a student of Terasaka in Osaka University and
the first Japanese to publish research on Knot Theory \cite{Kin-4}.

The first subject he studied was the fixed point property.  These
studies were concluded with the publication of \cite{Kin-1,Kin-3}.
After that, he became interested in topological transformations, 
which were, probably, the main achievements of Terasaka.
In April 1953 T.Homma came to Osaka University as a Yukawa
Fellow (compare entry on Homma). He stayed one year (till March 1954),
 and Terasaka and Kinoshita studied topological transformations with him.  
Homma and Kinoshita wrote three joint papers on this subject  
\cite{H-K-1,H-K-2,H-K-3}.
 After finishing their third paper they realized that it would
be difficult to advance the research further without studying 
3-dimensional topology, and they heard of the progress made by 
Moise, Bing and Harrold\footnote{Later he wrote a joint paper with 
Kinoshita \cite{H-K}.} in USA.
H.Noguchi, who joined the group, had been interested in the Poincar\'e 
Conjecture.  He strongly suggested to study knot theory and 
to contact Professor R.H.Fox at Princeton University. So, Kinoshita wrote a
letter to Fox asking for reprints of his papers\footnote{New journals
from abroad were available only in the library of American Cultural Center
\cite{Mim}, and there were no copying machines then.}.
In 1955, with Fox's papers in hand, they (Terasaka, Kinoshita and Yajima)
started a seminar on knot theory. 
In 1956, H.Noguchi came to Osaka University as a Yukawa fellow 
(April-September). The seminar grew and new participants included Noguchi,
Miss. Y.Hashizume, F.Hosokawa and J.Tao.
At first, they were interested in trivial knots (the Smith Conjecture, 
motivated by \cite{M-S}, and Dehn's Lemma).  
However, Terasaka advised Kinoshita to study
non-trivial knots (so to treat knots as the main object of study). 
Kinoshita received a doctorate in mathematics in 1958 for results
published in \cite{Kin-7} and announced in \cite{Kin-5}, and in 
1958-59 was a Lecturer at Osaka University.
In September of 1959 he moved to Princeton. For two years ((1959-61) he was a
member of the Institute for Advanced Study and later (1961-62) he was
a Research Associate at Princeton University. In these three years 
he collaborated closely with Professor Fox.
For two years (1962-64) he was Assistant Professor at
University of Saskatchewan (Canada). 
In 1964-68 he was an Associate Professor at Florida State University, and
1968 he became a Professor there and stayed at FSU till his 
retirement in December, 1984.
In January 1985 he moved back to Japan
to become a professor at Kwansei Gakuin University, where his
friend Yajima (also student of Terasaka) was a professor.
He retired from Kwansei Gakuin University in April of 1994.
Now Kinoshita lives in the suburbs of Atlanta (USA) with his wife Michiko.
They have 4 children: 3 girls and a boy.

His students include:\\
Yaichi Shinohara (b. 1942), PhD in 1969 at Florida State University
(now Dean at Kwansei Gakuin University).

At Kwansei Gakuin University, Kinoshita's students include:\\
Toshio Harikae, Hisanori Naka, Masahiro Nakao and Yasushi Yonezawa.

\item[(N)] Hiroshi Noguchi, born December 26, 1925 in Tokyo.\\
He spent two years attending night courses at Science University of 
Tokyo\footnote{At that time it was called Tokyo Science School, 
one of the professional training schools.  After attending such a school 
for about 3 years (possibly 2 during the war) 
one was able to apply to a degree awarding University.}
(working at daytime) during wartimes.
In 1945: ``At the age of nineteen I was a freshman in a college and
decided to study Topology as my major research field "\cite{Su-2}.
He graduated in 1948 from T\^{o}hoku University under 
Prof. Tadao Tannaka supervision.
He read Alexandrov-Hopf {\it Topology}. 
He was an instructor at T\^{o}hoku University and later moved 
to Waseda University (1952) and became an assistant professor. 
He started teaching from Analytic Geometry and Projective
Geometry. In 1953 he gave a lecture on topological spaces. 
In 1956 he got a Yukawa fellowship 
and went to work with Terasaka in Osaka.
In Osaka he met Kinoshita once again\footnote{
Noguchi met Kinoshita before, in 1954, at the Meeting of 
the Japanese Mathematical Society in Osaka. 
Then Kinoshita introduced Homma to Noguchi. Also in 1954
Homma and Kinoshita published their first joint paper \cite{H-K-1},
and it was reviewed by Fox. The first joint paper by Homma and Terasaka 
was published in 1953 \cite{Hom-Ter}.}. 
Mr and  Ms Tao (Ms Tao formerly Hashizume, known for prime
decomposition of links 1958) and Hosokawa. 
Noguchi spent one year in Osaka and went back for one year to Waseda.
There was new (1957) exchange program between Waseda and 
University of Michigan. Noguchi was one of the first people to use it. 
He spent a year at University of Michigan at Ann Arbor.
In January of 1958, in Michigan, he studied under E.Moise 
(there were there also R.Bott and H.Samelson).
In July of 1958 he went to Princeton (for one week). There he met 
Fox\footnote{When Noguchi met Fox he was asked about Murasugi. The reason was
that Murasugi had sent his papers, \cite{Mur-1}, to Fox. Fox told
Noguchi that the results were obtained before by his student R.H.Crowell.
Later Fox wrote joint review of the papers by Murasugi and Crowell,
and advised Murasugi to come to Princeton.}, 
 Papakyriakopoulos, Stallings and L.P.Neuwirth. 
Before that July he submitted a paper to the Transactions of the American
Mathematical Society but the referee (Moise) suggested that
he rewrite and split the paper. One part went to Osaka Journal 
and one to Annals of Mathematics 
- the best Mathematical Journal in the world  
(``Smoothing of combinatorial $n$-manifolds in $n+1$ space").
He got PhD for these two papers in 1960 at Osaka University.
In Michigan he met M{\"u}ller and in 1963/4 they published 
joint papers in computer science. 

Noguchi was a  Senior Foreign Scientist Fellow in Mathematics visiting
the University of Illinois in 1966-1967. 
After 1970 Noguchi became interested in Mathematical Education and published  
a high school text in 1975. He is the president of Japanese branch of 
the Mathematical Olympiad.

Noguchi dedicates his book \cite{No-3} to Tadao Tannaka 
(his supervisor at T\^{o}hoku University) and Komatu.

In Section 4 we list the students of Noguchi.

\item[(H)] Tatsuo Homma, born January 1926.\\
He graduated from Tokyo Institute of Technology\footnote{ Then called 
Tokyo Higher Technical School, one of the professional training schools.}
 under Minagawa\footnote{Takizo Minagawa was a student of
Terasaka specializing in differential geometry; compare \cite{Hom-Mi}.}. 
He got a Yukawa fellowship in 1953\footnote{Hideki Yukawa (Ogawa), 1907-81,
 was the first Japanese Nobel prize winner, 1949,
for his work on elementary particles (prediction of mezons). Hideki
Yukawa contributed a portion of the
prize to Osaka University, where he once held a position.  With this
contribution, the university established the Yukawa Fellowship.
 Ogawa assumed his wife's family name upon his marriage (1932) \cite{Yuk}.}
and went to work with Terasaka in Osaka (compare entry on Kinoshita). 
He attacked Dehn's Lemma in 1950's, first proved it partially
\cite{Hom-1}, and later in full generality (for $S^3$) \cite{Hom-2}, 
independently from Papakyriakopoulos. 
Homma was a professor at Yokohama City University, then at
Tokyo Institute of Technology, and later at Aoyama Gakuin University (Tokyo).
He was a member at the Institute for Advanced Study in Princeton (1960-62).
He was also a visiting associate professor at Florida State University
(1965), where he wrote Lecture Notes \cite{Hom-3}
and at Michigan State University (1967).
He gave an invited (1 hour) talk at the AMS meeting in Arkansas.
In 1973 he was on the organizing committee of the International 
Conference on Manifolds and related topic (the major topological conference). 
He was the chief organizer (from Japanese side)
of a joint Japan-US conference held at Hawaii in 1982.


His students include Yasuyuki Tsukui, Mitsuyuki Ochiai, Teruo Nagase 
and Seiya Negami.

\item[(M)] Kunio Murasugi, born March 25, 1929 in Tokyo.\\
Murasugi attended the Tokyo Higher Normal School from 1945, 
graduating in 1949 (compare the footnote 5). Then he entered the Tokyo Bunrika 
University\footnote{Tokyo Bunrika University was absorbed into the Tokyo 
University of Education in 1949
(Nobel Prize winner in physics, 1965, Sin-Itiro Tomonaga was
its president, 1956-1962).  In 1973 it changed the
official name to the University of Tsukuba.}
and graduated in 1952 with B.Sc. under Takeuti\footnote{Gaisi 
Takeuti was born in 1926 and graduated
from Tokyo University under
the supervision of Sh\^{o}kichi Iyanaga who specialized
in number theory, but had among his students also Kodaira. Takeuti's
main specialization was mathematical logic, however he was interested
in topology and about 1950 he worked on
homotopy theory. He proved a nice result, accepted for publication in
J.of Math. Society of Japan, but to his disappointment
the result just has been proven by G.W.Whitehead.
Then he moved his attention to knot theory.} supervision.
As a university student Murasugi wanted to study algebraic topology,
in particular, homotopy theory, and to determine the m-th homotopy group of
an n-sphere. This was a fashionable problem in early 1950's.
But the famous J.P. Serre's paper (about fiber spaces)
 appeared in Annals of Mathematics in 1951 \cite{Ser}.
Murasugi was discouraged by this paper and almost had lost his interest
in the problem. Takeuti was then interested in algebraic topology, so
Murasugi joined his weekly seminar on topology (he was the only student
in the seminar). After Murasugi graduated, Takeuti suggested to him that
in order to study homotopy theory he should join some active group
in the topic. He should go to Osaka and join the Toda group.
Otherwise, Takeuti continued, Murasugi would better study ``knot theory".
This was the first time Murasugi have heard  on knot theory \cite{Mur-4}.
In 1952 (or 1953), just after Murasugi's graduation in 1952, Takeuti
gave a semester long course in knot theory. It was based mostly on
Reidemeister's book \cite{Re} but also on his own thoughts
(he was then unaware of the work of Fox and his students\footnote{
Later, in Princeton, Takeuti played {\it Go} with Fox but
never mentioned his knot theory course \cite{Tak}.}).
Soon, Takeuti went back to his research in logic,
but Murasugi continued working in knot theory.
First he wanted to prove the asphericity of the knot complement, using 
group theoretical methods.
He found that the alternating knot group had some regularity,
and so he tried to solve the above problem for alternating knots.
However the problem was solved  by Aumann \cite{Au} for
alternating knots, and later by Papakyriakopoulos, for general 
knots \cite{Pap}. Nevertheless, Murasugi was able to find interesting
facts on alternating knots resulting in his
first three publications\footnote{Murasugi sent a reprint of
his first publication to Fox, who noticed that the similar result
was obtained by his student Crowell. The work of Crowell was then
unknown in Japan, even to Kinoshita who was a referee of Murasugi's 
paper.} in 1958 \cite{Mur-1,Mur-2},
which later become basis of his doctorate\footnote{The person 
with several published and well-received papers could
submit them as his doctoral thesis to his (or any other) university.
There was no need for a formal supervisor or formal course credits.
However Murasugi always considered Takeuti as his supervisor.}.

Murasugi became an assistant at Housei University in 1954,  
Lecturer\footnote{In the Japanese system there are three 
tenured positions: Lecturer, Assistant Professor and  Professor.  
They roughly correspond to American: Assistant Professor, Associate 
Professor and (full) Professor.}
in 1956, and Assistant Professor in 1960. Around 1958, Homma
 was working at the Tokyo Institute of Technology
and Noguchi was at Waseda University, also in Tokyo.
They used to have a weekly seminar at Homma's office at TIT.  
Since Noguchi lived near Murasugi's home, (most likely) he introduced 
Murasugi to Homma some time in 1958, after Murasugi published
his first paper. Murasugi first met Terasaka and Kinoshita at 
the Topology Conference in early summer of 1958. 
Kinoshita left Japan to Princeton in the summer of 1959 and Homma
went to Princeton in 1960. 

Murasugi met Moise and Bing in Japan\footnote{Murasugi remembers that
he, Homma and Noguchi met Moise and Bing at the Tokyo (Haneda)
airport.} before he left in 1960 to America, 
but he met Fox for the first time only in Princeton in 1961. 
One should add that J.P.Mayberry\footnote{PhD thesis, Princeton Univ., 
1955. Now he is a Professor Emeritus of 
Mathematics at Brock University in St. Catherines, Ontario, Canada.} 
often visited
Japan in late 1950's being involved with military duties in Korea. 
In particular he  gave a colloquium at Waseda University (invited 
by Noguchi), attended, in particular, by Murasugi\footnote{Later they
wrote a joint paper \cite{M-M}.}, Tao and Hashizume. 
In 1957 Coxeter and Moser published a book on infinite group 
theory \cite{C-M}.  Murasugi studied the book and found it very 
useful in his research in knot groups. His idea was to first study
infinite group theory with Coxeter\footnote{Harold Scott MacDonald
Coxeter (b.1907).}
 in Toronto and then study knot theory with Fox in Princeton.
He obtained a graduate scholarship as an MA student 
at the University of Toronto.
In 1960 Murasugi went to Canada\footnote{Murasugi arrived at Toronto
September 15, after 2 days flight from Japan. He didn't
take a ship as stated in \cite{No-3}.}.
In December, 1960, he was told that his
doctoral thesis submitted to Tokyo University of Education before he flew 
to Toronto was accepted and in March, 1961, Murasugi received the doctoral
degree from the Tokyo U.E. 
After getting a doctoral degree from Japan, Murasugi's status in Toronto 
changed.  He was no longer a graduate student but became, in April of 1961, 
a research assistant of William Thomas Tutte (b.1917), a specialist in 
graph theory. 
In November, 1960 (2 months after Murasugi left Japan), R.Fox
wrote him a letter in which he suggested Murasugi to apply for a visiting 
membership at the Institute for Advanced Study in Princeton 
(and enclosed an application form)\footnote{Next year, in November
 1961, Murasugi drove to Princeton to meet Fox for the first time 
and the whole his family stayed at Kinoshita's home for one week. He also
met Crowell there.}.
He consulted Professor
Gilbert de Beauregard Robinson (1906--1992) (representation theory) and
was strongly suggested to stay in Toronto. Murasugi agreed with him,
and it was one of the reasons that after 2 years at Princeton he
returned to Toronto not to Japan.
Murasugi was a research assistant of Tutte from April, 1961 to June, 1962.
In August,1962, he (and his family) went to Princeton University as a
research associate, and stayed there until August 1964. He got
an Assistant Professor position (in Toronto), in July 1964, became 
an associate professor in July 1966, and professor in July 1969.
Toronto became the place for young Japanese mathematicians\footnote{And
not only Japanese: I also was a postdoc of Murasugi 
in Toronto 13 years ago.}  to learn
Knot Theory and to establish relations with American topologists. 

Murasugi's students include:\\
Richard T.Hartley University of Toronto 1976 \cite{Har},\\
Toru Maeda\footnote{He is now working at Kansai University in Osaka.}, 
 University of Toronto, 1983 \cite{Ma},\\
Bohdan I. Kurpita, University of Toronto 1992,\\
Peter Clifford Hill, University of Toronto, 1998,\\
John Mighton, University of Toronto, 1999.

In July of 1999, there was a conference on Knot Theory in Toronto 
to celebrate Prof. Kunio Murasugi's 70th birthday \cite{Mur-3}.

\item[(H)] 
Fujitsugu Hosokawa (b.1930), graduated at Osaka University
 under the supervision of Terasaka.  
He moved to Kobe University.  Hosokawa is well known for
his paper on ``Hosokawa polynomial" (\cite{Hos}, 1957).
His father was also a mathematician. He wrote the first 
textbook about geometry (university Textbook) [T\^{o}hoku University]. 
Terasaka read his book.
Hosokawa senior moved to Hiroshima University and was killed in Hiroshima 
by the atomic bomb (at the time he was then conducting a seminar and
 some of the young students survived when the building collapsed). 
F.Hosokawa is a creator of the Kansai school of knot theory.
His students at Kobe University include:\\
Akio Kawauchi (master degree student\footnote{
Kawauchi was a student of
Taniguchi and Terasaka at Sophia University (Terasaka introduced Kawauchi
to knot theory and Taniguchi supervised seminar based on the Spanier's book). 
He did Master degree at Kobe University under Hosokawa and Suzuki 
supervision. He did his
doctoral degree under Junzo Tao at Osaka City University.
He spent 2 years at the Institute
for Advanced Study in Princeton (1978-1980).
He is a professor at Osaka City University. 
In 1983 he co-organized the KOOK (acronym for the main
participating Universities in Kansai district: Kobe University,
Osaka University, Osaka City University, and Kwansei
Gakuin University), geometric topology monthly seminar, which
educated several young mathematicians (see entry on J.Tao).
His students include: 
Makoto Sakuma (Doctorate), Hitoshi Murakami (Master, Dr), Adrian 
Pizer (Ms, Dr), (She) Masako Kobayashi (Dr), Seiichi Kamada (Ms, Dr), 
Yasuyuki Miyazawa (Dr), Teruhisa Kadokami (Dr), 
Shin Satoh (Ms, Dr), and Ikuo Tayama (Dr).}), 
Masakazu Teragaito (1986), Yoko Nakagawa (master)\footnote{
She did her PhD. under Fox at Princeton University and was his last
doctorate student. She is now working at Yamaguchi Women's
University.}, 
Tetsuo Shibuya, Kazuo Yokoyama, Yasutaka Nakanishi (he was one of 
the first doctoral students at Kobe University), Makoto Sakuma, 
Taizo Kanenobu, Yoshiaki Uchida, Yoshihiko Marumoto and Kanji Morimoto
 (doctorate).

\item[(Has)] Youko Hashizume, later Ms. Tao (b.1930)\footnote{
Several papers in Kobe Journal of Mathematics 6(1), 1989 
are dedicated to her on her 60th birthday.}, graduated at 
Osaka University under the supervision of Terasaka. She got a position at 
Kansai University. Her best known result concerns the uniqueness of 
the connected sum decomposition of links \cite{Has}.\footnote{I had been 
generalizing her result on free action on $S^1\times S^2$, \cite{Tao,Pr}.} 
 
\item[(Ta)] Junzo Tao (b.1929), 
graduated at Osaka University under the supervision of
Terasaka. His major field is differential topology.
In 1950's Tao and H.Yamasuge shared office at Osaka University.
Yamasuge tried to solve the Poincar\'e Conjecture, and presented 
his prove of 5-dimensional case \cite{Yam},
but died young, in November 1960.
Tamura, Kudo, Tao and Saito published his remaining papers.
Smale sent a dedicated letter about Yamasuge's work
to Tao and Yoshihiro Saito (1930-1997). K.Kobayashi, Kawauchi, 
Nakanishi, and Hitoshi Murakami were Tao's doctoral students. 
Tao was an early
member of Osaka City knot theory seminar (from 1956), he never 
wrote a paper on knots but his broad knowledge of topology
(including algebraic and geometric topology) was appreciated by
participants of the seminar.

KOOK Seminar  was established in 1983 by Tao and Kawauchi of Osaka
City University  with strong support  from F.Hosokawa of
Kobe University and S.Kinoshita (from January 1985 at Kwansei 
Gakuin University).
The impulse of this creation was the fact that a small conference 
hall of Osaka City University was opened in a convenient place in
Osaka City. It was felt that it would be fruitful to have 
regular meetings for 
scholars in knot theory (mostly students and friend of H.Terasaka).
KOOK Seminar was organized by geometric topologists from 
Kobe University, Osaka University, Osaka City University and
Kwansei Gakuin University, and was attended by many topologists,
mostly from Kwansai area. The naming of KOOK Seminar is due to Adrian
Pizer who was a student of Kawauchi \cite{Ka}. From 1994, geometric
topologists of Nara Women's University joined the seminar as
organizers (and since then it is called often the N-KOOK Seminar).

\item[(Ya)] Takeshi Yajima (1914-1998) was a pioneer of 
4-dimensional knot theory \cite{Yaj}.
Yajima was a student of Terasaka at Osaka University. Later he was
a professor at Kwansei Gakuin University.
His students at Kwansei Gakuin University include:\\
Kouhei Asano(MS, PhD),
Toru Maeda(MS),
Katsuyuki Yoshikawa(MS 1980, PhD 1985).
\item[(Yan)] Takaaki Yanagawa (b.1935), \cite{Yan}. He was 
the last student of Terasaka at Osaka University. He joined Osaka 
University knot theory seminar in April, 1959 (compare entry on Kinoshita). 
He worked at Kobe University and
now he is a professor at Kansai University in Osaka (invited by 
Hashizume (Tao), to fill her position after she retired).
\end{enumerate}

\section{Third generation}
Students of Hiroshi Noguchi:
\begin{enumerate}

\item[(Ku)]
Keiko Kudo (later Kamae); born 1939.
She studied at Waseda University under the supervision of Noguchi,
and left (even before finishing her master program) to become
an Instructor at Kyushu University, invited by J.Tao. She went to the USA 
and was an instructor at the University of Illinois. Later she became
an assistant professor at Nihon University. She wrote a joint paper
with Noguchi in 1963 \cite{Ku-No}.

\item[(Ko)] 
Kazuaki Kobayashi (b. December 17, 1940), graduated from Waseda University
in March 1964.
Master degree from Waseda 1966 under Noguchi supervision. 
Doctorate: \ 
Osaka City University under the supervision of  J.Tao, 1972; \cite{Ko}. 
He went in April 1966 to Kobe University (Math. Dept. in Science Faculty) 
in 1968 moved to Math. Dept.  in Art Faculty, and to Hokkaido in 1971.
He is now a professor at the Tokyo Women
Christian University.
His students include:\\
Masaharu Kouno, Ken Sakai, Yoshiyuki Yokota, Yuko Yoshimatsu, Takako Kodate, 
 Chikako Toba, Kumiko Endo, 
Makiko Ishiwata, Noriko Imi, Kaneko Masuda, and Shizuka Mitaki.

\item[(S)]
Shin'ichi Suzuki (b. April 20, 1941), graduated from Waseda University
in March 1965.
MS, Waseda Univ. Mar. 1967 (with Noguchi). Ph.D. Waseda Univ. March 1974. 
Terasaka invited Suzuki to Sophia University as an instructor in 1967.
Suzuki moved to Kobe University after 2-3 years. Later he created a school
of Knot Theory at the Waseda University in Tokyo.

There is the conference this year (October 2001) to celebrate 
60th birthday of Kobayashi and Suzuki.

Suzuki's students include:\\
In Kobe:  
Akio Kawauchi (MS), Shuji Yamada\footnote{He did his master and doctoral 
degree at Osaka University under supervision of Minoru Nakaoka.},
 Yasutaka Nakanishi, Kanji Morimoto (master and doctorate), 
Makoto Sakuma, and Taizo Kanenobu.\\ 
In Waseda: Teruhiko Soma (PhD), Yoshiyuki Ohyama (PhD), Masao Hara (PhD),
Kazuo Yokoyama (MS at Kobe, doctorate at Waseda)  
(She) Miki Shimabara (now Miyauchi) (master),
  Kouki Taniyama (doctorate), (She) Yasuko  Ogushi (now Ohyama) (MS)
  Yoshiyuki Yokota (doctorate), Satoshi Yamashita (doctorate)
  Akira Yasuhara (doctorate), Toshiki Endo (MS), Keigo Makino (MS),  
  (She) Tomoe Motohashi (doctorate), (She) Miyuki Okamoto (doctorate) 
  [graduate of Tsuda College], Eishin Kawamoto (MS), 
  Makoto Ozawa (doctorate), and Tatsuya Tsukamoto (doctorate)\footnote{He 
also has been my PhD student at GWU, graduating in May 2000.}.

\item[(Ka)] 
Mitsuyoshi Kato (b. October 10, 1942), graduated Waseda Univ. Mar. 1964.
PhD. 1970\footnote{All graduate work of Kato was done at Waseda University; 
his doctoral theses were submitted to Tokyo University  and doctorate was
awarded
it in 1970.}. He was a student of Noguchi at Waseda University.
He visited the Institute for Advanced Study in Princeton (1968-70).
He is now a Professor at Kyushu University.
Kato's students include:\\
 Sadayoshi Kojima (Master March 1978)\footnote{ 
Kojima went to study for PhD at Columbia University (New York),
were he spent two years 1979-1981, and wrote his Ph.D. thesis under the
supervision of John Morgan. Kojima is now a professor at
Tokyo Institute of Technology. His students include: 
Shigeru Kitakubo (Ph.D 1992, Graph Theory),
        Yosuke Miyamoto (Master 1990, Hyperbolic Geometry),
        Yasushi Yamashita (Master 1991, Ph.D 1996, Hyperbolic Geometry),
        Mitsuhiko Takasawa (Ph.D 2000, Surface Automorphism), and
        Kazuhiro Ichihara (Ph.D 2000, 3-Manifolds).}, 
Toshitake Kohno (master 1980)\footnote{Kohno graduated in 1979 at
Tokyo University. He went to Nagoya, later 
to Kyushu University, and finally to Tokyo University. Toshie Takata
 was his student.},
Masayuki Yamasaki (Master 1978)\footnote{Yamasaki did his PhD in 1982 at
the Virginia Polytechnic Institute and State University under the
supervision of Frank Quinn.}, 
Kimihiko Motegi, (She) Haruko Nishi, and Yoshihisa Sato.

\item[(Fuk)] 
Masako Fukuda (b. March 27, 1942), graduated Waseda Univ. Mar. 1964.
PhD. 1974. Her husband Takuo Fukuda wrote with Noguchi a textbook on 
elementary catastrophy theory, 1976.

\item[(Fuj)]
Kiichi Fujino (b. Feb. 15, 1931), graduated Waseda Univ. Mar. 1955.
PhD. 1972. Computer Science.

\item[(Hat)]
Mitsuhiro Hattori (b. Feb 18, 1941), graduated Waseda Univ. Mar. 1964.
PhD. 1974. Computer Science (asynchronous circuits).

\item[(Hon)]
Masaru Honda (b. Dec 15, 1941), graduated Waseda Univ. Mar. 1965.
Applied statistics.

\item[(Ik)]
Hiroshi Ikeda (b. April 18, 1944), graduated Waseda Univ. Mar. 1965.
PhD. 1974. Combinatorial topology and graph theory.

\item[(It)] 
Ryuuichi Ito (b. April 18, 1944), graduated Waseda Univ. Mar. 1967.
PhD. 1983.

\item[(Kas)]
Takumi Kasai (b. June 3, 1946), graduated Waseda Univ. Mar. 1969.
PhD. 1977. Computational Complexity.

\item[(Ma)]
Seishi Makino (b. Nov. 21, 1961), graduated Waseda Univ. Mar. 1987.
Computer Science.

\item[(Mo)]
Etsuro Moriya (b. Jan. 7, 1947), graduated Waseda Univ. Mar. 1970.
PhD. 1976. Computer Science.

\item[(Nakai)] Nakai (b. Feb. 24, 1957), graduated Waseda Univ. Mar. 1979.
PhD. 1985 (On topological types of polynomial map germs, Kyoto Univ).

\item[(Nakam)] 
Tsuyoshi Nakamura (b. March 3, 1947),  graduated Waseda Univ. Mar. 1969.
PhD. 1982 (Waseda Univ.), DMS. 1983 (Nagasaki Univ.).
Statistical Sciences.

\item[(Nar)]
Hiroshi Narushima (b. Sept. 30, 1943), graduated Waseda Univ. Mar. 1967,
D.Sc. Dec. 1977 (Waseda Univ.). Discrete Mathematics,  Computer Science.

\item[(Nishin)]
Tetsuro Nishino (b. Feb. 3, 1959) graduated Waseda Univ. Mar. 1982,
D.Sc. Dec. 1991.  Computer Science.

\item[(Nishim)] 
Takashi Nishimura (b. Feb. 5, 1958), graduated Waseda Univ. Mar. 1983,
D.Sc. Dec. 1988.  Singularity theory.

\item[(Sa)]
Ken Sawada (b. Feb. 28, 1953), graduated Waseda Univ. Mar. 1975,
D.Sc. Dec. 1981, Waseda Univ. (Extended orbits of diffeomorphisms
and Omega-explosions). 

\item[(So)]
Teruhiko Soma (b. June 26, 1955), graduated Waseda Univ. Mar. 1980.
 
\item[(Su)]
Shao-Chin Sung (b. Sep. 8, 1969), graduated Waseda Univ. Mar. 1993,
Computer Science.

\item[(Ta)]
Sei'ichi Tani (b. June 5, 1963), graduated Waseda Univ. Mar. 1987.
Ph.D. 1996. Computer Science.

\item[(To)]
Yoshio Togawa (b. Jan. 1953), graduated Waseda Univ. Mar. 1975,
D.Sc. Mar. 1978, Waseda Univ. Dynamical Systems.

\item[(Tsu)]
Kensei Tsuchida, (b. Mar. 29, 1958), graduated Waseda Univ. Mar. 1982.
D.S. 1994. Computer Science.

\item[(Yak)]
Takeo Yaku (b. Oct. 21, 1947), graduated Jiyu Gakuen College Mar. 1970.
M.Sc. Waseda Univ. Mar. 1972. Ph.D. Waseda Univ. Dec. 1977. Computer Science.

\item[(Yam)]
Makoto Yamamoto (b. Jan. 26, 1953), graduated Waseda Univ. Mar. 1976.
Doctor of Science 1982. Topological graph theory.

\end{enumerate}

\section{Other mathematicians whose work is related to knot theory}
Of course not all mathematicians interested in knot theory are
students or grand students of Terasaka. We give below three examples
of a world famous mathematician whose work is related to knot theory:
\begin{enumerate}
\item[(Mats)] Yukio Matsumoto (b. 1944) graduated from Tokyo University in 
1967 under the supervision of Tamura\footnote{Ichiro 
Tamura, well known specialists in Differential
Structures and  Characteristic Classes, was a student of 
Sh\^{o}kichi Iyanaga (see footnote 5), who was in turn 
a student of Teiji Takagi (1875--1960), specialist in class field theory, 
one of the first great Japanese mathematicians. 
Takagi studied in Germany (in Gottingen,
1899-1901, under Hilbert and Klein) sent there by the Japanese Government
as was the norm for young able people 100 years ago. His teacher was
 Rikitaro Fujisawa  (1861-1933) who represented
Japan at the ICM in  1900, \cite{Iy}.}.
Matsumoto visited the Institute for Advanced Study (in 76-78)
where he met Jos\'e Montesinos (leading Spanish mathematician) and they
started a very fruitful collaboration. Matsumoto is now 
the president of Mathematical Society of Japan.
Matsumoto's students include:\\
Ken'ichi Kuga, Zjunici Iwase, Katura Miyazaki, Masahico Saito 
(master\footnote{Saito and Miyazaki obtained their PhD
at University of Texas at Austin under the supervision of C.Gordon.}),
Ken'ichi Ohshika, Osamu Saeki (1992, PhD, Tokyo University), 
Kiyoshi Ohba, Kazumasa Ikeda, Nariya Kawazumi, Tomotada Ohtsuki, 
 Kazunori Kikuchi, Kazuhiko Kiyono, Chuichiro Hayashi, 
Toru Ikeda, Goro Fujita, Kazuo Habiro, Yuichi Yamada,  
 Akiko Shima, and Shigeru Takamura.

\item [(Fu)]
Shinji Fukuhara (now at Tsuda College) was born in 1945 
and graduated in March 1968 from Tokyo University, 
under the supervision of Ichiro Tamura. 
His students include:\\
Noriko Maruyama (Dr. Sci., 1994), Jinko Kanno, 
Haruko Miyazawa (formerly Aida, Dr. Sci., 2000),
Miyuki Okamoto (doctorate under S.Suzuki at Waseda University),
 Naoko Kamada (formerly Ishii), Kazuko Onda.

\item [(Matu)]
Takao Matumoto (now Hiroshima University). Born 1946, graduated in 1968
(Tokyo University) under the supervision of Ichiro Tamura. 
Then his advisor was Tokushi Nakamura (M1) and he obtained his master degree
at Tokyo University under the supervision of Akio Hattori (in 1970).
He did his doctorate in France (Orsay) under the supervision of Laurent 
C.Siebenmann. Mituhiro Sekine was Matumoto's student working in Knot Theory.  

\end{enumerate}

\section{Names of some Japanese topologists born before 1946}
\begin{enumerate}
\item[] Kuniz\^o Yoneyama (1877-1968),\  Takeo Wada (1882-1944),
\item[] Keitar\^{o} Haratomi (1895-1968),\  K\^oshiro Nakamura (1901-1985), 
\item[] Kinjir\^o Kunugui (1903-1975), Hidetaka Terasaka (1904-1996), 
\item[] Hiroshi Okamura (1905-1948), Sh\^{o}kichi Iyanaga (1906--)
\item[] Motokiti Kond\^{o} (1906-1980), \ Atuo Komatu (1909-1995),
\item[] Takeshi Inagaki (1911-1989),\ Shizuo Kakutani (1911--),
\item[] Kiyoshi Aoki (1913--), \      Makoto Abe (1914-1945)
\item[] Takeshi Yajima (1914-1998), \ Kiiti Morita (1915-1995),
\item[] Kunihiko Kodaira (1915-1997)\ Ryozi Sakata (Shizuma) (1916-1985),
\item[] Tatsuzi Kudo (1919--), \      Hiroshi Uehara (1923--),
\item[] Hidekazu Wada (1924--), \     Shin'ichi Kinoshita (1925--),
\item[] Hiroshi Noguchi (1925--), \   Jun-iti Nagata (1925--),
\item[] Kei\^o Nagami (1925--), \     Minoru Nakaoka (1925--),
\item[] Nobuo Shimada (1925--), \     Tatsuo Homma (1926--),
\item[] Ichiro Tamura (1926-1991), \  Hiroshi Yamasuge (1926-1960),
\item[] Katsuhiko Mizuno (1926--), \  Ichiro Yokota (1926--),
\item[] Hirosi Toda (1928--), \       Masahiro Sugawara (1928--),
\item[] Ken-ichi Shiraiwa (1928--), \ Kunio Murasugi (1929--),
\item[] Junzo Tao (1929--), \         Tsuneyo Yamanoshita (1929--),
\item[] Yoshihiro Saito (1930-1997),\ Fujitsugu Hosokawa (1930--),
\item[] Youko Hashizume (Tao) (1930--), \ Nobuo Yoneda (1930-1996),
\item[] Sh\={o}r\={o} Araki (1930--), \ Tokusi Nakamura (1930--),
\item[] K\^{o}zi Shiga (1930--), \      Haruo Suzuki (1931--),
\item[] Masahisa Adachi (1931-1993), \  Yasutoshi Nomura (1932--),
\item[] Seiya Sasao (1933--), \         Kunio \^{O}guchi (1933--),
\item[] Takaaki Yanagawa (1935--), \  Hiromichi Matsunaga (1935--),
\item[] Yoshihiro Shikata (1936--), \ Teiichi Kobayashi (1936--),
\item[] Fuichi Uchida (1938--), \     Mamoru Mimura (1938--),
\item[] Keiko Kudo (Kamae) (1939--) \ Kazuaki Kobayashi (1940--),
\item[] Shin'ichi Suzuki (1941--), \  Mitsuyoshi Kato (1942--),
\item[] Masako Fukuda (1942--), \     Katsuo Kawakubo (1942-1999),
\item[] Goro Nishida (1942--), \      Akihiro Tsuchiya (1942--),
\item[] Yaichi Shinohara (1942--), \  Hajime Sato (1944--),
\item[] Hiroshi Ikeda (1944--), \ Ryuuichi Ito (1944--),
\item[] Yukio Matsumoto (1944--),  Shinji Fukuhara  (1945 --) 
\end{enumerate}

\section{Ralph H.Fox}\label{7}
Ralph Hartzler Fox was born March 24, 1913.
A native of Morrisville, Pa., he attended Swarthmore
College for two years while studying piano 
at the Leefson Conservatory of Music in Philadelphia \cite{C-F,Fo-8}.
He received his master's degree from the Johns Hopkins University.
He obtained his Ph.D. from the Princeton University in 1939 under
the supervision of Solomon Lefschetz (1884--1972).
Fox was married, when he was still a student, to Cynthia Atkinson. 
They had one son, Robin.
After receiving his Princeton
doctorate, he spent the following year at
Institute for Advanced Study in Princeton.
He taught at the
University of Illinois and Syracuse University before
returning to join the Princeton University faculty in 1945
and staying there until his death in 1973.

He devoted
most of his career to mathematical topology, and in particular
to knot theory: the study of different ways of placing closed
curves or loops in three-dimensional space.

At the International Congress of Mathematicians, Cambridge, Mass., 1950,
Fox gave a talk on ``Recent development of knot theory at Princeton" 
\cite{Fo-7}.
He gave a series of lectures at the Instituto de Matem\'aticas
de la Universidad Nacional Aut\'onoma de M\'exico in the
summer of 1951. He lectured to the American Mathematical Society (1949),
to the Summer Seminar of the Canadian Mathematical Society (1953), and at the
Universities of Delft and Stockholm, while on a Fulbright 
grant (1952) \cite{Fo-3}.
His influential book, {\it Introduction to Knot Theory} \cite{C-F},
was based upon lectures given at Haverford College under the
Philips Lecture Program (spring 1956).
The last papers of Fox were published in Fundamenta Mathematicae
and Osaka Journal of Mathematics \cite{Fo-5,Fo-6}.

Fox career is summarized nicely in the book dedicated to his memory
\cite{Neu}: ``The influence of a great teacher and a superb 
mathematician is measured by his published works, the published
works of his students, and perhaps foremost, the mathematical environment
he fostered and helped to maintain."

One of Fox's interests was the ancient Japanese board
game of Go (which he learned as a graduate student).
He represented the United States in the first
international Go tournament, held in Tokyo in 1963, and
later received the fourth Dan degree conferred by the
international Go organization in Tokyo.

Fox was an active Quaker.

He died December 23, 1973 in the University
of Pennsylvania Graduate Hospital, where he had undergone open-hearth surgery.
He was survived by his widow, the former
Cynthia Atkinson\footnote{She died recently.}, 
and their son, Robin H. of Minneapolis \cite{Fo-8,Fo-10}.

Fox's Ph.D. students include:\\
A.L.Blakers (the first student of Fox, 1947), R.C.Blanchfield, 
W.A.Blankinship, E.J.Brody, R.H.Crowell, D.M.Dahm, C.H.Giffen, H.R.Gluck, 
D.L.Goldsmith, F.J.Gonz\'ales-Acuna, H.W.Kuhn, R.H.Kyle, 
S.J.Lomonaco, Jr., S.B.Mauer, Y.Nakagawa (the last student of Fox, 1973), 
B.C.Mazur, J.W.Milnor,  L.P.Neuwirth, N.F.Smythe, J.R.Stallings, Jr., 
G.Torres-Diaz, W.C.Whitten and E.F.Whittlesey. 
A.L.Anger, M.Artin\footnote{Michael Artin is a son of Emil Artin,
and he did his Ph.D under the supervision of Oscar Zariski.}, D.A.Gay and
Kenneth A.Perko, Jr. wrote their Princeton Senior Thesis under 
the supervision of Fox.

\section{Acknowledgments}\label{8}
I am grateful to Prof. Kazuaki Kobayashi for his invaluable help in 
preparing this essay. Most of the work was done when I was visiting 
TWCU in the summer of 1999. I would like to thank my Japanese hosts
for their hospitality.

\centerline{Department of Mathematics,}
\centerline{ The George Washington University,}
\centerline{Washington, DC 20052}
\centerline{e-mail: przytyck@gwu.edu}
\end{document}